\documentclass[12pt,reqno]{amsart}
\usepackage{times}
\usepackage{amssymb}
\usepackage{fullpage} 
\usepackage{setspace}
\usepackage{mathtools}
\usepackage{comment}
\usepackage{color}
\usepackage{amsmath,amsfonts,amsthm,amssymb,amscd,url}
\usepackage{graphicx}
\usepackage{mathrsfs}
\usepackage{todonotes}
 \usepackage{dsfont}
\usepackage{bbold}

\usepackage[parfill]{parskip} 

\usepackage[toc,page]{appendix}

\usepackage[utf8]{inputenc}
\usepackage{fullpage,amssymb,amsmath}
\usepackage{xcolor}
\usepackage{halloweenmath}
\usepackage{url}
\usepackage{hyperref}
\newtheorem{theorem}{Theorem}

\title{Moments of $L$-functions Problem List}
\author{
Alia Hamieh, Habiba Kadiri, Greg Martin, Nathan Ng}
\date{last updated \today}

\begin{document}

\maketitle

\section*{Introduction}
This is an ongoing list of problems that has resulted from the PIMS (Pacific Institute of Mathematical Sciences) Collaborative Research Group
\href{https://www.pims.math.ca/collaborative-research-groups/lfunctions}{
 $L$-functions in Analytic Number Theory: 2022-2025}. 
 The main focuses of this CRG include moments of $L$-functions and automorphic forms, explicit results in analytic number theory, and comparative prime number theory.
% A full list of activities are available at the 
%\href{https://sites.google.com/view/crgl-functions/crg-news}{CRG website}. 
This list began at the workshop \href{https://sites.google.com/view/crgl-functions/workshop-moments-of-l-functions#}{\it Moments of $L$-functions} which took place at the University of Northern British Columbia from 
July 25, 2022 to July 29, 2022.  Other problems were added at the  CRG launch meeting
\href{https://www.birs.ca/events/2022/2-day-workshops/22w2003}{\it L-function in Analytic Number Theory}
which took place at the Banff International Research Station from Nov. 18, 2022 to Nov. 20, 2022. 
The main goal of this list is to stimulate research in this field. 
%The CRG also has an active research \href{https://sites.google.com/view/crgl-functions/crg-weekly-seminar?authuser=0}{seminar}. 
%Videos for the seminar are available on PIMS video sharing site \href{https://www.mathtube.org/CRG/167}{mathtube}.  \\

The focus of this problem list is on moments of $L$-functions and related topics.  The study of Moments of $L$-functions has a long history.
The classic problems in the field  are related to the $2k$-th moments of the Riemann zeta function, 
\begin{equation}
 \label{IkT}
I_k(T) = \int_{0}^{T} |\zeta(\tfrac{1}{2}+it)|^{2k} \, dt. 
\end{equation}
 A foundational book in the field is Titchmarsh's book  {\it The Theory of the Riemann zeta function.}  This book outlines topics on the size, the distribution, and the moments  of the Riemann zeta function.  The main conjecture in the field is the Riemann Hypothesis, the assertion that all non-trivial zeros 
of the Riemann zeta function line on the one-half line. Another main conjecture is the  Lindel\"{o}f Hypothesis, the statement that
\begin{equation}
\forall \varepsilon >0, \,  \, |\zeta(\tfrac{1}{2}+it)| \ll (1+|t|)^{\varepsilon}.
\end{equation}
  Hardy and Littlewood devised an approach to resolving this conjecture which involved establishing, for every $\varepsilon >0$, there exists $C(\varepsilon) >0$ 
 \begin{equation}
 I_k(T) \le C(\varepsilon)  T^{1+\varepsilon}. 
 \end{equation}
 Indeed, in 1918, they succeeded in evaluating the second moment and demonstrated 
 \begin{equation}
  \label{secondmoment}
   I_1(T) \sim  T (\log T). 
 \end{equation}
 In 1926, Ingham evaluated the fourth moment and showed that 
 \begin{equation}
  \label{fourthmoment}
    I_2(T) \sim \frac{T}{2 \pi^2} (\log T)^4. 
 \end{equation}
In Titchmarsh's book \cite[Chapter 7]{Ti}, it is shown  that 
\[
   \int_{0}^{\infty} e^{-t/T} |\zeta(\tfrac{1}{2}+it)|^{2k} \, dt  \gg  T (\log T)^{k^2}. 
\]
This is a smoothing of $I_k(T)$ and strongly suggests that the correct size is 
\begin{equation}
  \label{IkTasymptotic}
   I_k(T) \sim C_k T (\log T)^{k^2}
\end{equation}
for a certain constant $C_k$.   Despite the success in determining the second moment \eqref{secondmoment} and the fourth moment  \eqref{fourthmoment}, 
$I_k(T)$ has not been asymptotically evaluated for any other moment. 

In the 1960's  moments of Dirichlet $L$-functions were studied.
For instance, in Montgomery's monograph \cite[Theorem 10.1, p.75]{M} a  bound for
\begin{equation}
 \sum_{\chi} \int_{-T}^{T} |L(\tfrac{1}{2}+it,\chi)|^{2k} \, dt  
\end{equation}
was established in the case $k=2$.  The related moments \footnote{To be precise, Bombieri studies smoothed versions of these moments.}
\begin{equation}
 \sum_{q \le Q}   \sum_{\chi} \int_{-T}^{T}  |L(\tfrac{1}{2}+it,\chi)|^{2k} \, dt 
\end{equation}
are studied in Bombieri's book \cite{Bo} on the large sieve inequality. 
These results were proven because of their connection to zero-density results and to primes in arithmetic progressions.  
Dirichlet $L$-functions were also studied at their central point.  This is related to Chowla's conjecture that quadratic Dirichlet $L$-functions 
do not vanish at $s=1/2$.  In an attempt to establish the non-vanishing of Dirichlet $L$-functions at the central point, researchers examined mean values of the shape
\begin{equation}
  \label{quad}
  \sum_{|d| \le X}  L(1/2,\chi_d)^k
\end{equation}
where $\chi_d$ ranges through real Dirichlet characters.  
They also studied the moments 
\begin{equation}
  \label{charq}
  \sum_{\chi \text{ mod } q} L(1/2,\chi)^k
\end{equation}
where $\chi$ ranges over all Dirichlet character modulo $q$. 

In the 1980's researchers began to study moments of  higher degree $L$-functions, including those attached to automorphic forms.
It was discovered that the vanishing or non-vanishing of such $L$-functions have arithmetic consequences.  For instance, Kolyvagin \cite{K}
showed that the finiteness of the Tate-Shaferevich group was the consequence of demonstrating that certain elliptic curve $L$-functions
had simple zeros.  This led to independent research by Murty-Murty \cite{MM} and Bump-Friedberg-Hoffstein \cite{BFH} where they evaluated 
the moments
\begin{equation}
  \label{elliptic}
  \sum_{|d| \le X} L'(1/2, E \otimes \chi_d)^k
\end{equation}
where $L(s,E)$ is the normalized $L$-function associated to an an elliptic curve $L$-function.  The method used in \cite{MM} is the technique of Approximate Functional
Equations and the method in \cite{BFH} is the technique of Multiple Dirichlet Series.  These are two of the main techniques in studying
moments of $L$-functions.   Another application from studying modular $L$-functions may be found in the work of Ellenberg \cite{E}. He showed 
that the non-vanishing of certain modular $L$-functions could be used in showing that the Fermat equation $A^4+B^2=C^p$ (where $p$ is prime),  has no non-trivial solutions. 

In the 1990's some new ideas entered this field of study.   In 1998 Keating and Snaith \cite{KS} developed a random matrix model for moments of the Riemann zeta function. Using their model they came up with a precise conjecture for the constants $C_k$ in \eqref{IkTasymptotic}, for all $k >0$.    They showed that 
\begin{equation}
 \label{Ckgkak}
 C_k =  \frac{G(1+k)^2}{G(1+2k)} a_k
\end{equation}
where $G$ is the Barnes $G$-function and 
\begin{equation}
 \label{ak}
 a_k=
 \prod_{p} \Big( 
  (1-p^{-1})^{k^2} \sum_{m=0}^{\infty} \Big( \frac{\Gamma(m+k)}{m! \Gamma(k))} \Big)^2  p^{-m}
  \Big).
\end{equation}
For a very long time researchers had been unable to determine the value of $C_k$.  Around the same time Katz and Sarnak \cite{KS2} introduced the notion of symmetry type in families of  $L$-functions.    In their book they showed that the distribution of the zeros of certain function field $L$-functions follow 
a law first discovered by Montgomery and then verified numerically by Odlyzko.  In  \cite{KS2}  the authors consider moments of the shape
\begin{equation}
  \label{general}
  \sum_{\substack{f \in \mathcal{F} \\ c(f) \le X}} L(1/2,f)
\end{equation}
where $\mathcal{F}$ is a ``family" of $L$-functions and $c(f)$ denotes the ``conductor" of $f$.  Note that the moments of 
\eqref{quad}, \eqref{charq}, and \eqref{elliptic} are all special cases of \eqref{general}.  In \cite{KS2} a conjecture for the asymptotic size of the moments
\eqref{general} was formulated also based on the random matrix model ideas introduced in \cite{KS}.  
The conjectural formulae for these moments depend on the symmetry type of each family $\mathcal{F}$. 
Before this article, researchers did not have a good idea what the correct answer was for these moments.  The conjecture of Keating-Snaith in \cite{KS2} also agreed with all proven asymptotic  formulae at the time. These conjectures  generated a significant
amount of research and a flurry of activity.   The random matrix model method is useful as it gives a  straightforward way to 
conjecturally study problems on moments of $L$-functions.   One defect is that it does not see the ``arithmetic" in moment calculations. 
Later is was discovered that analytic number theory methods could lead to the
same conjectures as in \cite{KS}, \cite{KS2}. For instance, the article \cite{CFKRS} uses the approximate functional equation (``the recipe")
to provide a full main term conjecture for any reasonable family of  moments including \eqref{IkT} and \eqref{general} 
for $k$ a natural number. 
In \cite{DGH}, the authors derive the leading term in the moment conjectures using the method of Multiple Dirichlet Series.
The connection between the Approximate Functional Equation method
and the Multiple Dirichlet Series method is not well understood. However, there is currently work in progress by Baluyot and \u{C}ech
studying this (see also \cite{Ce}).   Ideas on the approximate functional equation side have been developed by Conrey and Gonek \cite{CG}.
They modelled moments of the zeta function by mean values of long Dirichlet polynomials. In particular, they gave a heuristic calculation showing
how derive the sixth and eighth moments of the Riemann zeta function \cite{CG} from these mean values. 
Their heuristic calculations have 
 been made rigorous
in the articles \cite{Ng}, \cite{NSW}, and \cite{HN}. 
Conrey and Keating recently wrote a series of article \cite{CK1}, \cite{CK2}, \cite{CK3}, \cite{CK4}, and \cite{CK5} which attempt to generalize
the method of Conrey-Gonek to higher moments.   In \cite{CK4} and \cite{CK5} they introduce a new type of additive divisor sum.  
This is an active field of study and  researchers are trying to understand this new approach and these new divisor sums (see \cite{BC} for example).

  Another major theme in this field is the connection to the spectral theory of automorphic forms.
Ideas of this nature originated in work of Kuznetsov.  It was systematically developed by Deshouillers and Iwaniec in \cite{DI}. 
In this article they outlined many applications of the spectral theory to automorphic forms to analytic number theory. 
In particular, they derived sharp bounds for additive divisor sums in \cite{DI2} and applied their methods to moments of the Riemann zeta function
\cite{DI3}, \cite{DI4}.  Motohashi applied the spectral theory to additive divisor sums in \cite{M} and then produced an exact formula for the smoothed
 fourth moment of the Riemann zeta function in the monograph \cite{Mo}.  His surprising formula relates the fourth moment to a third moment of 
 Maass form $L$-functions.  Ideas from Motohashi's work are actively being studied now other researchers have tried to find analogous formulae 
 for other moments (for instance, see \cite{BHKM} and \cite{HK2}).  Kevin Kwan's recent series of article \cite{Kw1}, \cite{Kw2}, and \cite{Kw3} provide a useful background
 for this topic. 
 
 This introduction is just a broad overview of this field and is not intended to be a comprehensive survey of the field. 
There are now a number of excellent surveys on moments of $L$-functions including \cite{C} and  \cite{S2} which the reader may consult.

If someone has solved or made progress on a problem or if someone has a problem to contribute to this list 
 please contact {\texttt{alia.hamieh@unbc.ca}} or {\texttt{nathan.ng@uleth.ca}}. 
 If you would like to contribute a problem, please send a clearly written problem with any  references in latex code. 
 If you have solved one of the problems 
 from this list please acknowledge and cite this article. 

\section*{Notation}
\begin{itemize}
    \item For $k \ge 0$, define 
    \[
     I_k(T) = \int_{0}^{T} |\zeta(\tfrac{1}{2}+i t)|^{2k} \, dt.
    \]
 \item For $k >0$, $d_k(n)$ denotes the $k$-th divisor function, defined by 
 \[
   \zeta(s)^k = \sum_{n=1}^{\infty} \frac{d_k(n)}{n^s}.
 \]
 \item For $k,\ell >0$, define for $x >0$ and $h \in \mathbb{N}$
 \[
   D_{k,\ell}(x,h) = \sum_{n \le x} d_k(n) d_{\ell}(n+h)
 \]
\end{itemize}

\section*{Problems}

\subsection*{Moments of the Riemann zeta function}
\begin{enumerate}
\item Give an explicit version of Harper's \cite{Ha} upper bound for $I_k(T)$, assuming the Riemann hypothesis.  That is, for $k >0$ determine $C_k$ such that 
\begin{equation}
 \label{ub}
  I_k(T) \le C_k T (\log T)^{k^2}
\end{equation}
for $T$ sufficiently large. 
 Note that this is related to the following conjecture of Soundararajan \cite{S2}
for $k \ge 2$ and $T \ge 10$,
\[
   I_k(T) \le  T (\log T)^{k^2}.
\]
(Nathan Ng, Lethbridge) \\
\begin{itemize}
\item S. Inoue \cite{In} has established a form of \eqref{ub}.
\item We have also been informed a graduate student has worked extensively on this problem and nearly has a result. 
\end{itemize}

\item Smoothed fourth moment.  Let $\omega(t)$ be a smooth compactly supported function supported in $[c_1T, c_2T]$ where $0 < c_1 < c_2$. 
Show that 
\[
  \int_{-\infty}^{\infty} \omega(t) |\zeta(\tfrac{1}{2}+i t)|^{4} \, dt = O(T^{\frac{1}{2} + \theta+\varepsilon  })
\]
where 
$\theta$ is upper 
bound in Ramanujan's conjecture. 
See \cite[Corollary 2]{Ber} where an analogous result was proven  for modular $L$-functions.

\item In 1989 Zavorotnyi proved that 
\[
  I_2(T) = T P_4 (\log T) + O( T^{\frac{2}{3}+\varepsilon}),
\]
where $P_4$ is a specific degree $4$ polynomial. 
Improve the bound in the error term.  \\
(Nathan Ng, Lethbridge)

\item Conrey {\it el al} made the conjecture 
\[
  I_k(T) = TP_{k^2} (\log T) + O(T^{\frac{1}{2} + \varepsilon}) 
\]
where $P_{k^2}$ is a certain degree $k^2$ polynomial.  Let 
\[
 E_k(T) = I_k(T) -TP_{k^2} (\log T) 
\]
denote the error term.  Since the publication of \cite{CFKRS} there has been some debate on what should be the correct
conjectural size of 
$E_k(T) $. Note that it is widely believed that $E_1(T) = O(T^{\frac{1}{4}+\varepsilon})$ and 
$E_2(T) = O(T^{\frac{1}{2}+\varepsilon})$ are best possible.  There are omega results for both of these cases.  It would be desirable to come up with a conjecture  for optimal
constants $\theta_k$ such that 
\[
  E_k(T) = O(T^{\theta_k+\varepsilon})
\] 
holds for $k \in\mathbb{N}$. N. Ng believes that $\theta_1=\frac{1}{4}$ and $\theta_2=\frac{1}{2}$ are the correct answers for 
$k=1,2$. What about $k \ge 3$?  Both Ivic and Motohashi suggested $\theta_3 = \frac{3}{4}$.  Can the Multiple Dirichlet Series method be used to derive
an appropriate conjecture? (Work in progress of Baluyot and Cech is studying this.)  \\
(Nathan Ng, Lethbridge)

\item Improve Chandee's \cite[Theorem 1.2]{Ch}.
\begin{theorem}
Let $k_i$ be positive real numbers. Let $\alpha_i=\alpha_i(T)$ be real-valued functions of $T$ such that $\alpha_i =o(T)$. Assume that $\lim_{T\rightarrow \infty} \alpha_i\log T$ and $\lim_{T\rightarrow \infty} (\alpha_i-\alpha_j)\log T$ exist or equal $\pm \infty$. Assume that for $i\neq j$, $\alpha_i\neq \alpha_j$ and  $\alpha_i-\alpha_j =O(1)$. 
Then the Riemann Hypothesis implies that for $T$ sufficiently large, one has
$$
  M_{{\bf k}}(T,  {\bf \alpha})\ll_{{\bf k},\varepsilon} T(\log T)^{k_1^2 + \cdots k_m^2 +\varepsilon} \prod_{i<j} \left( \min\left\{\frac{1}{|\alpha_i-\alpha_j|}, \log T \right\}\right)^{2k_ik_j}
$$
where
\[
   M_{{\bf k}}(T,  {\bf \alpha})
   = \int_{0}^{T} |\zeta(\tfrac{1}{2} + (t+\alpha_1) )|^{2 k_1} \cdots |\zeta( \tfrac{1}{2} + (t+\alpha_m)|^{2 k_m}  dt.
\]
\end{theorem}
Use Harper's method \cite{Ha} to remove the $\varepsilon$. Note that the case of two shifts ($m=2$ and $k_1=k_2$) has been dealt with in \cite{NSW}. \\
(N. Ng, Lethbridge)
\begin{itemize}
\item Michael Curran  \cite{Cu} has solved this problem.
\end{itemize}

\item Soundararajan  \cite{S} proved the following lower bound for $I_3(T)$:
\[
  I_3(T) \ge \frac{20.26 a_3}{9!} (1+o(1)) T (\log T)^9,
\]
for $T$ sufficiently large.  Improve the value of 20.26. The best possible result would be 42. 
(Likewise improve the lower bounds for $I_k(T)$ for $k >3$.)  Note that the works \cite{BBLR}, \cite{BCR}, and \cite{HY} are now available and could be useful.  \\
(Nathan Ng, Lethbridge)
\begin{itemize}
\item We have been informed a graduate student has worked extensively on this problem and nearly has a result. 
\end{itemize}

\item Is it possible to prove an omega theorem for $I_3(T)$ assuming the ternary additive divisor conjecture?
Prove that there exists $\Theta >0$, such that 
\[
  I_3(T) = T \mathcal{P}_9(\log T) + \Omega(T^{\Theta})
\]
where $\mathcal{P}_9$ is a certain degree nine polynomial. \\
(Nathan Ng, Lethbridge)

\item Ng, Shen, and Wong \cite{NSW} showed that the Riemann hypothesis and a smoothed quaternary additive divisor conjecture imply 
the conjectured asymptotic 
\[
  I_4(T) \sim C_4 T (\log T)^{16}
\]
for a certain positive constant $C_4$.  Under the same assumptions, obtain the full main term, up to an error term of size $o(T)$ or a power savings error term.  Obtain a shifted version of this result.   Under the same assumptions, asymptotically evaluate
\[
  \int_{-\infty}^{\infty} \omega(t) \Big( \prod_{j=1}^{4}  \zeta(\tfrac{1}{2}+a_j+it) \zeta(\tfrac{1}{2}+b_j-it)  \Big) \, dt
\]
where $a_j, b_j$, $j=1, \ldots, 4$ are complex numbers satisfying $|a_j|, |b_j| \ll \frac{1}{\log T}$ and $\omega(t)$ is a smooth, compactly supported function in $[c_1T, c_2T]$ with $0 < c_1 < c_2$. 

\item In \cite{CFKRS} a conjectural method known as the `recipe' is introduced.  The `recipe' was used to produce  (conjecturally) a full asymptotic formula for $I_k(T)$.
Can the `recipe' be made rigorous in the case of the fourth moment?  That is, can we give a new proof of the asymptotic formula with power savings error term for $I_2(T)$
following the `recipe' method.  \\
(Nathan Ng, Lethbridge) \\

\item Can the conjectural lower order terms for $I_k(T)$ (due to CFKRS  \cite{CFKRS} and Diaconu-Goldfeld-Hoffstein \cite{DGH}) be recovered by refining the splitting conjecture in the Gonek-Hughes-Keating \cite{GHK} hybrid Euler-Hadamard product? \\
(Anurag Sahay, Rochester)

\item Prove third and fourth moment asymptotics with best possible error terms for $L(1/2,f)$ as $f$ runs over holomorphic Hecke eigenforms of weight $k_f \leq T$. It should be possible to show, using current technology, that
\[\sum_{k_f \leq T} \frac{L\left(\frac{1}{2},f\right)^3}{L(1,\operatorname{ad}f)} = T^2 P_3(\log T) + O_{\varepsilon}(T^{1+\varepsilon})\]
and
\[\sum_{k_f \leq T} \frac{L\left(\frac{1}{2},f\right)^4}{L(1,\operatorname{ad}f)} = T^2 P_6(\log T) + O_{\varepsilon}(T^{4/3+\varepsilon}),\]
where $P_j(x)$ denotes some polynomial of degree $j$ in $x$. For Hecke--Maass cusp forms, this is known for the fourth moment via work of Ivic, while Zhi Qi proved earlier this year the analogous result for the third moment.\\
(Peter Humphries, Virginia)

\item Heath-Brown proved the famous Weyl-strength twelfth moment bound
\[\int_{0}^{T} \left|\zeta\left(\frac{1}{2}+it\right)\right|^{12} \, dt \ll_{\varepsilon} T^{2 + \varepsilon}\]
for the Riemann zeta function. Jutila generalised this to prove an analogous result for the sixth moment in the $t$-aspect of $L(1/2 + it,f)$, where $f$ is a Hecke--Maass cusp form. Recently, Milicevic-White and Nunes have independently proven complementary results on the twelfth moment of Dirichlet $L$-functions in the level aspect. Can one also prove analogous bounds for the sixth moment of $L(1/2,f \otimes \chi)$?\\
(Peter Humphries, Virginia)

\item Jutila proved the Weyl-strength twelfth moment bound
\[\sum_{t_f \leq T} \frac{L\left(\frac{1}{2},f\right)^{12}}{L(1,\operatorname{ad}f)} \ll_{\varepsilon} T^{4 + \varepsilon},\]
where the sum is over Hecke--Maass cusp forms of level 1 and spectral parameter at most $T$. Can one prove a level-aspect analogue of this, where one sums over Hecke--Maass cusp forms of level $q$ and bounded spectral parameter?\\
(Peter Humphries, Virginia)

\item Good showed that given a Hecke--Maass cusp form of level 1 and fixed spectral parameter $t_f$, we have that
\[\int_{0}^{T} \left|L\left(\frac{1}{2}+it,f\right)\right|^2 \, dt = T P_1(\log T) + O_{t_f,\varepsilon}(T^{2/3 + \varepsilon}),\]
where $P_1(x)$ denotes some polynomial of degree 1 in x. Crucially, the dependence on $t_f$ in the error term is not made explicit. Can one make this dependence explicit? In particular, does this error term not depend
 on $t_f$ if, say, $T \geq t_f^{3/2 + \delta}$ for some $\delta > 0$? How about for holomorphic Hecke eigenforms of weight $k_f$?\\
(Peter Humphries, Virginia)

\item Use the recent work of  Humphries-Khan \cite{HK}  to compute weight aspect upper bounds for the 6th moment of $L(\frac12,f\otimes g)$ where $f$ and $g$ vary over the set of primitive forms of level 1 and weights $k$ and $\ell$ respectively, with $|k-\ell|=O(1)$ and $K<k+\ell\leq2K$.\\
 (Alia Hamieh, UNBC)

\item Remove the integral in $t$ from Chandee-Li's work and give an asymptotic formula for
$$
\frac{2}{\phi(q)} \sum_{\substack{\chi \mod{q}\\ \chi(-1)=(-1)^k}} \sum_{f\in\mathcal{H}_\chi}^h \left|L\left({\textstyle \frac{1}{2}},f\right)\right|^6
$$
See work of Chandee-Li-Matomaki-Radziwill on the sixth moment of Dirichlet $L$-functions.\\
(Joshua Stucky, University of Georgia)

\item Consider the family of Rankin-Selberg $L$-functions $L(\frac12,f\otimes g)$ where $g$ is fixed and $f$ varies over the set of primitive forms of weight $k$ and level $p^{\nu}$ as $\nu\to\infty$. Establish upper bounds for third and fourth moments of this family following the work in \cite{ChL1} and \cite{ChL2}.
Note that the 4th moment of the family $L(\frac12,f)$ with $f$ varying as above was studied by Balkanova \cite{balkanova}.  \\
 (Alia Hamieh, UNBC)

\item Matthew Young \cite{Y} proved that 
\begin{equation}
 \label{young}
  \sum_{0 < d < X}  L(\tfrac{1}{2}, \chi_{8d}) = X P(\log X) + O(X^{\frac{1}{2}+\varepsilon}). 
\end{equation}
where $d$ ranges through odd squarefree numbers and $\chi_{8d}$ denotes the quadratic character associated to the fundamental discriminant $8d$.  Alexandra Florea has suggested that 
\begin{equation}
 \label{floreaconj}
     \sum_{0 < d < X}  L(\tfrac{1}{2}, \chi_{8d}) = X P(\log X) + cX^{\frac{1}{4}} + O(X^{\vartheta})
\end{equation}
may be true for some $c \in \mathbb{R}$ and $\vartheta \in (0,\frac{1}{4})$.  The goal would be to improve \eqref{young} by reducing the exponent $\frac{1}{2}$. 
Note that Florea has proven a version of \eqref{floreaconj} in the function field case (see \cite{F}). \\
(Nathan Ng, Lethbridge)

\item Determine a polynomial $P$ of degree 10 such that 
\begin{equation}
 \label{young}
  \sum_{0 < d < X}  L(\tfrac{1}{2}, \chi_{8d})^4 = X P(\log X) + o(X) 
\end{equation}
where $d$ ranges through odd square-free numbers and $\chi_{8d}$ denotes the quadratic character associated to the fundamental discriminant $8d$.
Can the error term $o(X)$ be replaced by a power savings error term? Can the Multiple Dirichlet series method be used to obtain the first coefficients of $P$ (for instance see \cite{DPP}). \\
(Nathan Ng, Lethbridge)
\begin{itemize}
\item Shen and Stucky \cite{ShSt} have determined the first four coefficients of this polynomial. 
\end{itemize}

\item Compute moments for families of Artin-Schreier $L$-functions and improve results on 1-level zero density.\\
Notation: $p>2$ prime, $q$ power of $p$, $\mathcal L(u,f,\psi)$ is the Artin-Schreier $L$-function 
associated with rational function $f\in\mathbb F_q(x)$ and additive character $\psi$,
$$\mathcal{AS}^0_d=\{f\in\mathbb F_q[x]:\deg f=d\},\,(d,p)=1\}$$
$$\mathcal{AS}^{0,\mathrm{odd}}_d=\{f\in\mathbb F_q[x]:\deg f=d,\,f(x)=-f(-x)\},\,(d,2p)=1\}$$
\begin{multline*}\mathcal{AS}^{\mathrm{ord}}_d=\left\{f=h/g:h,g\in\mathbb F_q[x],(g,h)=1,g\mbox{ squarefree},\right.\\ \left. \deg f=\max(\deg h,\deg g)=d,\,\deg g\in\{d,d-1\}\right\}.\end{multline*}
$M_k=\langle \mathcal L(q^{1/2},f,\psi)^k\rangle$ is the $k$-th moment. All problems below concern the regime where $q$ is fixed and $d\to\infty$. 1-level density of zeros is defined w.r.t a Schwartz test function $\phi$.
\definecolor{dgreen}{RGB}{0,128,0}
Estimated difficulty scheme: {\color{dgreen}green = easy}, {\color{blue}blue = challenging}, {\color{red}red = hard}, black = $\bigskull$.
\\~\\
$\mathcal{AS}^{\mathrm{ord}}$: compute 1-level density for $\mathrm{supp}\,\hat\phi\subset{\color{red}(-1-\delta,1+\delta)\mbox{ for some }\delta>0}$ or $(-2-\delta,2+\delta),\mbox{ for some }\delta>0$.\\
Compute ${\color{dgreen}M_1},{\color{blue}M_2},{\color{red}M_3,M_4},M_5,M_6$.
\\~\\
$\mathcal{AS}^0$: compute 1-level density for $\mathrm{supp}\,\hat\phi\subset(-(2-2/p)-\delta,2-2/p+\delta)$.\\
Compute ${\color{dgreen}M_1,M_2,M_3},M_4$.
\\~\\
$\mathcal{AS}^{0,\mathrm{odd}}$: compute 1-level density for $\mathrm{supp}\,\hat\phi\subset(-(1-1/p)-\delta,1-1/p+\delta)$.\\
Compute ${\color{blue}M_1},M_2$.

(Alexei Entin, Tel Aviv University)

\item Investigate the dependence of zero-density  on moments of the zeta-function? More specifically, Ingham’s zero-density estimate makes use of the fourth moment $\int_{T}^{2T} |\zeta(1/2 + it)|^{4}\, dt$. Can we use other non-trivial moments (like, say, Heath-Brown’s estimate of the 12th moment) in such zero-density results, and can we generalise this to other L-functions?

(Timothy Trudgian, UNSW Canberra)

\section*{Bounds for $L$-functions}
\item Chandee and Soundararajan \cite{CS} proved that the Riemann Hypothesis implies
\[
  |\zeta(\tfrac{1}{2}+it)|
  \ll \exp \left( (C_0 +o(1)) \frac{\log t}{\log \log t} \right)
\]
with $C_0 = \frac{\log 2}{2}= 0.1505 \ldots $.  Reduce the value of $C_0$.  \\
(N. Ng, Lethbridge)

\item Let $S(t)=\frac{1}{\pi} \arg \zeta(\tfrac{1}{2}+it)$.   Carneiro, Chandee, and  Milinovich \cite{CCM} proved that the Riemann hypothesis 
implies
\[
  |S(t)| \le  (C_0+o(1)) \cdot \frac{ \log t}{\log \log t}
\]
with $C_0 =\frac{1}{4}$.  Reduce the value of $C_0$. \\
(Nathan Ng, Lethbridge)

\item Let $f$ be a Hecke--Maass cusp form of level 1. Prove the Weyl-strength subconvex bound $L(1/2+it,f) \ll_{\varepsilon} ((1 + |t_f - t|)(1 + |t_f + t|))^{1/6 + \varepsilon}$. This is known, via work of Jutila and Motohashi, for all $t \in \mathbb{R}$ away from the conductor-dropping regimes, namely when either $t_f - t$ or $t_f + t$ is $o(\max\{t_f,|t|\})$. Michel-Venkatesh prove a subconvex bound in the conductor-dropping regime, but the exponent is not explicit (and certainly not Weyl-strength).\\
(Peter Humphries, Virginia)

\section*{Additive divisor sums/shifted convolution sums}

\item  Topacogullari  \cite{T1}, \cite{T3} and Drappeau \cite{Dr} have proven asymptotic formulae for $D_{k,2}(x,h)$.
Prove an omega theorem for $D_{k,2}(x,h)-M_{k,2}(x,h)$ for $k \ge 3$ where $M_{k,2}(x,h)$ is the main term in the asymptotic formula. \\
(N. Ng, Lethbridge)

\item  Topacogullari \cite{T1}  proved that for $1 \le h \ll x^{\frac{2}{3}}$
\[
  D_{3,2}(x,h) = x P_3(\log x) + O(x^{\frac{8}{9}+\varepsilon}) 
\]
for a certain degree 3 polynomial where the coefficients of $P_3$ depend on $h$.   A result was also obtained for $D_{2,3}(x,h)$. 
\begin{itemize}
\item Improve the range of $h$ and the size of the error term. 
\item Compute exact simple formulae for the coefficients of $P_3$.  Numerically check the size of $| D_{3,2}(x,h) - x P_3(\log x)| $.
Is the error term of size $\sqrt{x}$?  See the recent related article of Nguyen \cite{Nguyen}.
\end{itemize}
(Nathan Ng, Lethbridge)

\item  Problems on additive divisor sums:
\begin{itemize}
\item Asymptotically evaluate 
\[
  \sum_{n \le x} d_3(n) d_3(n+1).
\]
(folklore)

\item Asymptotically evaluate
\[
  \sum_{n \le x} d(n) d(n+1) d(n+2).
\]
(folklore)
\item Let $\kappa \in (0,1)$. Asymptotically evaluate
\[
    \sum_{n \le x} d_{\kappa}(n) d(n+1) d(n+2).
\]
\item Let $\kappa \in (2,3)$. Asymptotically evaluate
\[
  \sum_{n \le x} d_{\kappa}(n) d_3(n+1). 
\]
\end{itemize}
\item Find a (conjectural) asymptotic formula for 
\[
  \sum_{n \le x} d_{k_1}(n+h_1) d_{k_2}(n+h_2) \cdots d_{k_r}(n+h_r)
\]
where $(k_1, \ldots, k_r) \in \mathbb{N}^r$ and $h_1, \ldots, h_r$ are distinct positive integers with $r \ge 3$. 
This would be analogous to the. conjectural asymptotic formula for 
\[
  \sum_{n \le x} \Lambda(n+h_1) \cdots \Lambda(n+h_r)
\]
that arises from the Hardy-Littlewood conjectures.  
Note that Ng-Thom and Tao have used different probabilistic arguments to derive conjectures in the case $r=2$.

\item Let $k  \in \mathbb{N}$ and $k \ge 3$.
For some small $\epsilon_0\geq0$ asymptotically evaluate
\[
  \int_{0}^{T} \sum_{m \le T^{2+\epsilon_0}} \frac{d_k(m)}{m^{\frac{1}{2}+it}}
  \sum_{n \le T^{2+\epsilon_0}} \frac{d(n)}{n^{\frac{1}{2}-it}}  \, dt.
\]
Note that for mean values of long Dirichlet polynomials we know how to evaluate these sums (conjecturally on the additive divisor conjecture) when 
the length of the polynomials is $N=T^{2-\epsilon}$ (see \cite{HN}).

\item Let $g$ be a primitive cusp form of weight $k$ for $\Gamma_{0}(N)$. Assume that $g$ is arithmetically normalized so that its first Fourier coefficient $a_g(1)=1$ and write $g(z)=\sum_{n=0}^{\infty}\lambda_{g}(n)n^{\frac{k-1}{2}}e(nz)$.  Let $X, Y, P\geq 1$, and let $F(x,y)$ be a smooth weight function supported  on $(X/4,4X]\times[Y/4,4Y]$ satisfying  \[x^iy^j\frac{\partial^{i}}{\partial{x^i}}\frac{\partial^{j}}{\partial v^j}f(u,v)\ll P^{i+j}\] for all $i,j\geq0$.
Consider the shifted convolution sum \[D_{g}(\ell_1,\ell_2;h)=\sum_{\ell_1m\pm \ell_2n=h}\lambda_g(n)\lambda_g(m)F(\ell_1m,\ell_2n),\] where $\ell_1$ and $\ell_2$ are coprime. Harcos \cite{Harcos} proved an asymptotic upper bound for $D_{g}(\ell_1,\ell_2,h)$ where the implicit constants depend on $g$. Blomer and Harcos  \cite{BH} established an asymptotic bound for an averaged version of this shifted convolution sum in which the implicit constants also depend on $g$.
\begin{itemize}
\item Determine the dependence of the explicit constant on $k$ and $N$. 
\item For a primitive cusp form $g$ of weight $k$ for the full modular group, establish an asymptotic upper bound that is uniform in the weight aspect. Such bounds are useful in moment computations for Rankin-Selberg $L$-functions associated with primitive forms of varying large weights.
\end{itemize}
 (A. Hamieh, UNBC)
% \item What is the conjecture for the error terms for $I_k(T)$?  CFKRS \cite{CFKRS} made the conjecture
% \[
%   I_k(T) = T \mathcal{P}_{k^2}(\log T)+ O(T^{\frac{1}{2}+\varepsion}).
% \]
% In some talks this has been revised to $O(T^{1-\delta_k+\varepsilon})$.  It seems likely that
% \[
%   \delta_1= \frac{1}{4} \text{ and }
%   \delta_2 = \frac{1}{2}.
% \]
% Does the Multiple Dirichlet Series method give a conjecture?

% \item Improve Drappeau's article.  What does GRH bound give?
% \item Improve error term in first moment of Quadratic Dirichlet $L$-functions. (Goldfeld-*)?
% \item Lower order terms in fourth moment of Quadratic Dirichlet $L$-functions.
\end{enumerate}

% \section{Introduction}

\noindent {\bf Acknowledgements}. We thank Alexei Entin, Peter Humphries, Anurag Sahay, Joshua Stucky, and Timothy Trudgian, 
for contributing problems to this list.

\end{document}